\documentclass[12pt]{article}

\usepackage{amssymb,amsmath,amsthm, amsfonts}

\def\a{{\mathrm{a}}}  \def\b{{\mathrm{b}}}  
    \def\f{{\mathrm{f}}}

\def\cA{{\mathcal{A}}}  
 \def\cE{{\mathcal{E}}} \def\cF{{\mathcal{F}}}
  
  \def\cL{{\mathcal{L}}}

 \def\cT{{\mathcal{T}}}

\let\le=\leqslant 
\let\ge=\geqslant 
\def\og{\leavevmode\raise.3ex\hbox{$\scriptscriptstyle\langle\!\langle$~}}
\def\fg{\leavevmode\raise.3ex\hbox{~$\!\scriptscriptstyle\,\rangle\!\rangle$}}
\def\essinf{\mathop{\rm ess\,inf}\limits}
\def\esssup{\mathop{\rm ess\,sup}\limits}

\def\sqr#1#2{{\vcenter{\vbox{\hrule height.#2pt
              \hbox{\vrule width.#2pt height#1pt \kern#1pt \vrule width.#2pt}
              \hrule height.#2pt}}}}
\def\signed #1{{\unskip\nobreak\hfil\penalty50
              \hskip2em\hbox{}\nobreak\hfil#1
              \parfillskip=0pt \finalhyphendemerits=0 \par}}
\def\endpf{\signed {$\sqr69$}}

\def\dbR{{\mathop{\rm l\negthinspace R}}}

\def\3n{\negthinspace \negthinspace \negthinspace }
\def\2n{\negthinspace \negthinspace }

\def\dbR{{\mathop{\rm l\negthinspace R}}}
\def\={\buildrel \triangle \over =}

\def\a{\alpha}
\def\b{\beta}

\def\f{\varphi}

\def\o{\omega}

\def\G{\Gamma}
\def\D{\Delta}

\def\F{\Phi}
\def\O{\Omega}

\def\cA{{\cal A}}

\def\cE{{\cal E}}
\def\cF{{\cal F}}

\def\cL{{\cal L}}

\def\cT{{\cal T}}

\def\ms{\medskip}

\def\qq{\qquad}

\def\esssup{\mathop{\rm esssup}}

\def\as{\hbox{\rm a.s.{ }}}

\def\esssup{\mathop{\hbox{\rm esssup}}}
\def\essinf{\mathop{\hbox{\rm essinf}}}

\def\|{\Big |}
\def\({\Big (}
\def\){\Big )}
\def\[{\Big[}
\def\]{\Big]}
\def\be{\begin{equation}}
\def\bel{\begin{equation}\label}
\def\ee{\end{equation}}
\def\bt{\begin{theorem}}
\def\et{\end{theorem}}
\def\bc{\begin{corollary}}
\def\ec{\end{corollary}}
\def\bl{\begin{lemma}}
\def\el{\end{lemma}}
\def\bp{\begin{proposition}}
\def\ep{\end{proposition}}
\def\br{\begin{remark}}
\def\er{\end{remark}}
\def\ba{\begin{array}}
\def\ea{\end{array}}
\def\ed{\end{document}}

\newtheorem{theorem}{Theorem}[section]
\newtheorem{lemma}[theorem]{Lemma}
\newtheorem{e-proposition}[theorem]{Proposition}
\newtheorem{corollary}[theorem]{Corollary}
\newtheorem{e-definition}[theorem]{Definition\rm}
\newtheorem{remark}{\it Remark\/}

\newtheorem{theoreme}{Th\'eor\`eme}[section]

\newtheorem{proposition}[theoreme]{Proposition}

\newtheorem{definition}[theoreme]{D\'efinition\rm}

\begin{document}

\title{Dual Representation as Stochastic Differential Games of Backward Stochastic Differential Equations and Dynamic Evaluations\thanks{This
work is partially supported by the NSFC under grants 10325101
(distinguished youth foundation) and 101310310 (key project), and
the Science Foundation of Chinese Ministry of Education under
grant 20030246004.}}

\author{ Shanjian Tang\thanks{Department of Finance and Control
Sciences, School of Mathematical Sciences, Fudan University,
Shanghai 200433, China, \& Key Laboratory of Mathematics  for
Nonlinear Sciences (Fudan University), Ministry of Education.
{\small\it E-mail:} {\small\tt sjtang@fudan.edu.cn}.\ms} }

\date{}

\maketitle \thispagestyle{empty}
\abstract{In this Note, assuming that the generator is uniform
Lipschitz in the unknown variables, we relate the solution of a
one dimensional backward stochastic differential equation with the
value process of a stochastic differential game. Under a
domination condition, an $\cF$-consistent evaluations is also
related to a stochastic differential game. This
 relation comes out of a min-max representation for uniform
Lipschitz functions as affine functions. The extension to
reflected backward stochastic differential equations is also
included.}

\section{Introduction}

Let $(\O,\cF, P)$ be a probability space, and $\{B_s;s\ge 0\}$ a
$d$-dimensional Brownian motion defined on $(\O,\cF, P)$. Let
$\cF_t$ be the  $\sigma$-algebra
 generated by  $\{B_s;0\le s\le t\}$ and the totality of $P$-null
 sets in  $\cF$, $L^2(\cF_t)$ the set of all $\cF_t$-measurable random variables $X$ such that $E|X|^2<\infty$,
 and $\cL^2_{\cF}(0,T)$ the set of $\cF_t$-adapted processes $\f$ such that $E\int_0^T|\f|^2\,
 dt<\infty$. Denote by $\cT_t$ the set of all
 $\cF_s$-stopping times taking values in $[t,T]$.

Consider the following  one dimensional backward stochastic
differential equation (BSDE): \bel{bsde} \left\{\ba{rcl}
dy_s&=&-f(s,y_s,z_s)\,
ds+\langle z_s, dB_s\rangle, \quad 0\le s\le T; \\
 y_T&=&\xi\in L^2(\cF_T).\ea\right.
\ee

It is known that when the generator is convex or concave with
respect to the unknown variables,  BSDE~(\ref{bsde}) is related
with a stochastic control problem.

More precisely, assume that $f$ is concave  in the last two
variables. Consider the Fenchel-Legender transformation: \be F(\o,
t,\b_1,\b_2):=\sup_{(y,z)}[f(\o,t,y,z)-\b_1 y-\langle
\b_2,z\rangle]\ee for any $ (\o,t,\b_1,\b_2)\in \O\times
[0,T]\times \dbR\times \dbR^d$. Define \be
D_t^F(\o)=\{(\b_1,\b_2)\in \dbR\times \dbR^d:
F(\o,t,\b_1,\b_2)<\infty\}.\ee Then the set $D_t^F$ is a.s.
bounded. It follows from well-known results (see, e.g.,
~\cite{ElKaPenQu}) that \be f(\o,t,y,z)=\inf_{(\b_1,\b_2)\in
D_t^F(\o)}[F(\o,t,\b_1,\b_2)+\b_1 y+\langle\b_2,z \rangle],\ee and
the infimum is achieved. Let us now denote by $\cA$ the set of
bounded progressively measurable $\dbR\times \dbR^d$ valued
processes $\{\b_1(t),\b_2(t)):0\le t\le T\}$ such that \be
E\int_0^TF(t,\b_1(t),\b_2(t))^2\, dt<\infty. \ee To each
$(\b_1,\b_2)\in \cA$, we associate the unique adapted solution
\break $(Y^{\b_1, \b_2},Z^{\b_1, \b_2})$ of BSDE~(\ref{bsde}) with
the coefficient $f$ being replaced with the affine one $f^{\b_1,
\b_2}(t,y,z):=F(t,\b_1(t),\b_2(t))+\b_1(t)y+\langle\b_2(t),z\rangle$.
In~\cite[pages 35--37]{ElKaPenQu}, the solution $y$ of
BSDE~(\ref{bsde}) is interpreted as the value process of a control
problem. That is, \be
y_t=\essinf_{(\b_1,\b_2)\in\cA}E[\F(t,\b_1,\b_2)|\cF_t]\ee where
\be
\F(t,\b_1,\b_2):=\D_{t,T}^{\b_1,\b_2}\xi+\int_t^T\D_{t,s}^{\b_1,\b_2}F(s,\b_1(s),\b_2(s))\,
ds\ee and for each $t\in [0,T]$, $\{\D_{t,s}^{\b_1,\b_2}: t\le
s\le T\}$ is the unique solution of the following stochastic
differential equation (SDE): \bel{SDE}
d\D_{t,s}=\D_{t,s}[\b_1(s)\, ds+\langle\b_2(s), dB_s\rangle],
\quad s\in [t,T]; \quad \D_{t,t}=1. \ee

The purpose of this Note is to obtain a similar dual
representation for the solution $y$ of BSDE~(\ref{bsde}) under the
Lipschitz assumption on the generator, instead of the convexity
assumption on the generator $f$.

Assume throughout the rest of the Note that there is a constant
$C>0$ such that \be \label{uniformLip} \left\{\ba {rl}\hbox{\rm
(B1)} & f(\cdot, y,z)\in \cL^2_{\cF}(0,T) \hbox{ \rm for any pair
} (y,z)\in
\dbR\times \dbR^d;\\
\hbox{\rm (B2)} & |f(t,y_1,z_1)-f(t,y_2,z_2)|\le
C(|y_1-y_2|+|z_1-z_2|) \\
&\hbox{\rm  for any } t\in [0,T]\hbox{ \rm and } (y_1,z_1),
(y_2,z_2)\in \dbR\times \dbR^d. \ea\right.\ee Then, for any $X\in
L^2(\cF_t)$, there is unique adapted solution $\{(Y_s,Z_s);0\le
s\le t\}$ of BSDE~(\ref{bsde}) with the terminal condition:
$Y_t=X$. Define $\cE^f_{s,t}[X]:=Y_s$ for any $s\in [0,t]$.

The rest of this Note is organized as follows. In section 2, we
give a Min-Max representation of a Lipschitz function in terms of
affine functions, which is the basis of the Note. In Section 3, we
present the dual formula for the solution of one dimensional
BSDE~(\ref{bsde}). In Section 4, the formula obtained in Section 3
 is applied to the dynamical evaluation and a dual formula is
 therefore derived  for an $\cF_t$-consistent evaluation. Finally in Section 5,
 a dual formula is also obtained for one dimensional reflected
 backward stochastic differential equations (RBSDEs)~(\ref{rbsde}).

\section{Min-max representation of a Lipschitz function as affine functions}\label{s2}

The following representation is due to Evans and
Souganidis~\cite[pages 786--787]{EvansSoug}.

\bl Let $f: [0,T]\times \Omega\times \dbR^n\to \dbR$ be a
Lipschitz function. That is, there is a constant $C>0$ such that
\be |f(t,x_1)-f(t,x_2)|\le C|x_1-x_2|, \quad \forall x_1, x_2\in
\dbR^n. \label{Lipschitz}\ee Then for each $t\in [0,T]$ and $x\in
\dbR^n$, \be f(t,x)=\max_{z\in \dbR^n}\min_{y\in {\overline
O}_n(0,1)}\{C\langle y, x\rangle+F(t,y,z)\}\ee
 where $F(t,y,z):=f(t,z)-C\langle y,z\rangle$ for $y,z\in \dbR^n$ and ${\overline O}_n(0,1)$ is the closed unit ball in $\dbR^n$.
 \el

{\bf Proof. } In view of the assumption~(\ref{Lipschitz}), we have
for any $x\in \dbR^n$ $$\ba{rl} f(t,x)=&\displaystyle\max_{z\in
\dbR^n}\{f(t,z)-C|x-z|\}\\
=&\displaystyle \max_{z\in \dbR^n}\min_{y\in {\overline
O}_n(0,1)}\{f(t,z)+C\langle y, x-z\rangle \}.\ea$$
\endpf

 \br
See Fleming~\cite[pages 996--1000]{Flem} or Evans~\cite{Evans} for
other, more complicated ways of writing a nonlinear function as
the max-min (or min-max) of affine mappings.
 \er



\section{Backward stochastic differential equations and related stochastic differential games}

Denote $(\cL^2_\cF(0,T))^{d+1}$ by $\cL^2_\cF(0,T;\dbR^{d+1})$,
and by $V_{d+1}$  the subset of $\cL^2_\cF(0,T;\dbR^{d+1})$ whose
element takes values in the closed unit ball ${\overline
O}_{d+1}(0,1)$.

Define the function $F: \O\times [0,T]\times \dbR^{d+1}\times
\dbR^{d+1}\to \dbR$ as follows: \bel{F} F(\o,s,\b_1,\b_2, \a_1,
\a_2)=f(\o, s,\a_1, \a_2)-C\b_1\a_1-C\langle \b_2, \a_2\rangle\ee
for any $(\o,s,\b_1,\b_2,\a_1,\a_2)\in \O\times [0,T]\times
\dbR^{d+1}\times \dbR^{d+1}$. Then, in view of Lemma 2.1, we have
for any $(\o,s,\b_1,\b_2,\a_1,\a_2)\in \O\times [0,T]\times
\dbR^{d+1}\times \dbR^{d+1}$, \be f(\o,t,y,z)=\max_{\a\in
\dbR^{d+1}}\min_{\b\in {\overline
O}_{d+1}(0,1)}[F(\o,t,\b,\a)+C\b_1y+C\langle \b_2, z\rangle].\ee

Given $\a\in \cL^2_\cF(0,T;\dbR^{d+1})$ and $\b\in V_{d+1}$,
consider the related BSDE: \bel{lbsde} \left\{\ba{rcl}
dY_s&=&-[C\b_1(s) Y_s+C\langle \b_2(s), Z_s\rangle\\
&& +F(s,\b_1(s),\b_2(s), \a_1(s),
\a_2(s))]\, ds+\langle Z_s, dB_s\rangle; \\
Y_T&=&\xi\in L^2(\cF_T).\ea\right. \ee
The solution is denoted by
$(Y^{\a,\b},Z^{\a,\b})$ when it is necessary to emphasize the
dependence on $(\a,\b)$ with $\a=(\a_1,\a_2)$ and
$\b=(\b_1,\b_2)$.

Introduce the following stochastic differential equation (SDE):
\bel{lSDE} d\G_{t,s}=\G_{t,s}[C\b_1(s)\, ds+C\langle\b_2(s),
dB_s\rangle], \quad s\in [t,T]; \quad \G_{t,t}=1. \ee Its solution
is denoted by $\G_{t,s}^\b, t\le s\le T$ to indicate the
dependence on $\b=(\b_1,\b_2)$.

 We have
 \bel{dual1}
 Y_t^{\a,\b}=E\left[\int_t^T\G_{t,s}^{\b}F(s, \b_1(s), \b_2(s), \a_1(s), \a_2(s))\,
 ds+\G_{t,T}^\b\xi\biggm |\, \cF_t\right]
 \ee for any $t\in [0,T]$.

\bt Assume that the function $f$ satisfies~(\ref{uniformLip}). Let
$(y,z)$ be the adapted solution of BSDE~(\ref{bsde}) and
$\{\G^{\b}_{t,s}; t\le s\le T\}$ the solution of SDE~(\ref{lSDE}).
Then we have for any $t\in [0,T]$,  \be
\ba{rcl}y_t&=&\displaystyle \esssup_{\a\in
\cL^2_\cF(0,T;\dbR^{d+1})} \essinf_{\b\in
V_{d+1}}E\biggl[\int_t^T\G_{t,s}^{\b}F(s, \b_1(s), \b_2(s),
\a_1(s), \a_2(s))\,
 ds\\
 &&\displaystyle\qq\qq\qq\qq\qq\qq\qq\qq+\G_{t,T}^\b\xi\biggm |\, \cF_t\biggr].\ea \ee \et

\section{An $\cF_t$-consistent evaluations and its dual representation as a stochastic differential game}

\begin{definition} A system of operators $\cE_{s,t}:L^2(\cF_t)\to L^2(\cF_s),
0\le s\le t\le T$ is called an $\cF_t$-consistent evaluation
defined on $[0,T]$ if it satisfies the following four properties:
for any $0\le s\le t\le T$ and any $X_1,X_2\in L^2(\cF_t)$,

(A1) $\cE_{s,t}[X_1]\ge \cE_{s,t}[X_2]$, \as, if $X_1\ge X_2$,
\as;

(A2) $\cE_{t,t}[X_1]=X_1$, \as;

(A3) $\cE_{r,s}[\cE_{s,t}[X_1]]=\cE_{r,t}[X_1]$, \as;

(A4) $\chi_A\cE_{s,t}[X_1]=\chi_A\cE_{s,t}[\chi_AX_1]$, \as for
any $A\in \cF_s$.

\end{definition}

In view of Peng~\cite[Corollary 4.2, page 588]{Peng1}, the
following is an immediate consequence of Theorem~\ref{fth1}.

 \bt Let $\{\cE_{s,t}\}_{0\le s\le
t\le T}$ denote an $\cF_t$-consistent evaluation  defined on
$[0,T]$. Assume that there is a function $g_\mu(t,y,z):=\mu
(|y|+|z|), (t,y,z)\in [0,T]\times \dbR\times \dbR^d$ for some
$\mu>0$ such that the $\cF_t$-consistent evaluation
$\{\cE_{s,t}\}_{0\le s\le t\le T}$ is dominated by
$\cE_{s,t}^{g_\mu}$ in the following sense: for any $s,t\in [0,T]$
such that $s\le t$ and for any $X_1,X_2\in \cL^2(\cF_t)$, we have
\be\cE_{s,t}[X_1]-\cE_{s,t}[X_2]\le \cE_{s,t}^{g_\mu}[X_1-X_2],
\quad \as.\ee Furthermore, assume that there is  $g_0\in
\cL^2_\cF(0,T)$ such that
$$\cE_{s,t}^{-g_\mu+g_0}[0]\le \cE_{s,t}[0]\le \cE_{s,t}^{g_\mu+g_0}[0].
$$
Then there is a function  $f:\Omega\times [0,T]\times
\dbR^{d+1}\to \dbR$ which satisfies~(\ref{uniformLip}), such that
$$
\cE_{s,t}[\xi]=\esssup_{\a\in \cL^2_\cF(0,T;\dbR^{d+1})}
\essinf_{\b\in V_{d+1}}
E\left[\G_{s,t}^{\b}\xi+\int_s^t\G_{s,r}^\b F(r,\a(r),\b(r))\,
dr\biggm|\, \cF_s\right].$$
 Here $\{\G^{\b}_{t,s}; t\le
s\le T\}$ is the solution of SDE~(\ref{lSDE}) and the function
$F:\Omega\times [0,T]\times \dbR^{d+1}\times \dbR^{d+1}\to \dbR$
is given by~(\ref{F}).

 \et

\section{Reflected backward stochastic differential equations and related mixed stochastic differential games}

We make the following assumption.

(B3) The obstacle $\{S_t, 0\le t\le T\}$ is a continuous
progressively measurable real-valued process satisfying \be
E\sup_{0\le t\le T}(S_t^+)^2<\infty, \quad S_T\le \xi, \as. \ee

Consider the following RBSDE: \be \label{rbsde} \left\{\ba{rcl}
dy_t&=&-f(t,y_t,z_t)\,
dt-da_t+\langle z_t, dB_t\rangle; \\
 y_T&=&\displaystyle \xi\in L^2(\cF_T); \quad y_t\ge S_t, \as
 \forall t\in [0,T];\quad
\int_0^T(y_t-S_t)da_t=0.\ea\right. \ee In view of \cite[Theorem
5.2, page 718]{ElKaPaPenQu}, it has a unique solution $(y,z,a)$.

Given $\a\in \cL^2_\cF(0,T;\dbR^{d+1})$ and $\b\in V_{d+1}$,
identically as in Section 3, consider the function $F$ given
by~(\ref{F}) and the related RBSDE: \bel{lrbsde} \left\{\ba{rcl}
dY_s&=&\displaystyle -[C\b_1(s) Y_s+C\langle \b_2(s), Z_s\rangle
\\
&&\displaystyle +F(s,\b_1(s),\b_2(s), \a_1(s),
\a_2(s))]\, ds-dA_s+\langle Z_s, dB_s\rangle; \\
Y_T&=&\displaystyle \xi\in L^2(\cF_T); \quad y_s\ge S_s, \as
 \forall s\in [0,T];\\
 &&\displaystyle \quad
\int_0^T(Y_s-S_s)dA_s=0.\ea\right. \ee
 The unique solution is
denoted by $(Y^{\a,\b},Z^{\a,\b},A^{\a,\b})$. We have for any
$t\in [0,T]$,
 \bel{dual3}\ba{rcl}
 Y_t^{\a,\b}&=&\displaystyle \esssup_{\tau\in  \cT_t}E\biggl[\int_t^\tau\G^\b_{t,s}F(s, \b_1(s), \b_2(s), \a_1(s), \a_2(s))\,
 ds\\
 &&\displaystyle \qq\qq\qq\qq+\G^\b_{t,\tau} S_\tau\chi_{\{\tau<T\}}+\G^\b_{t,\tau}\xi\chi_{\{\tau=T\}}\biggm |\cF_t\biggr].
\ea \ee

\bt  Assume that the function $f$ satisfies~(\ref{uniformLip}) and
the obstacle $\{S_t, 0\le t\le T\}$ satisfies assumption (B3). Let
$(y,z,a)$ be the adapted solution of RBSDE~(\ref{rbsde}) and
$\{\G^{\b}_{t,s}; t\le s\le T\}$ the solution of SDE~(\ref{lSDE}).
Then we have for any $t\in [0,T]$, \be \ba{c} \displaystyle
y_t=\esssup_{\a\in \cL^2_\cF(0,T;\dbR^{d+1}), \tau\in
\cT_t}\essinf_{\b\in V_{d+1}} E\biggl[\int_t^\tau\G_{t,s}^\b F(s,
\b_1(s), \b_2(s), \a_1(s), \a_2(s))\,
 ds\\
  \displaystyle+\G^\b_{t,\tau} S_\tau\chi_{\{\tau<T\}}+\G_{t,\tau}^\b\xi\chi_{\{\tau=T\}}\biggm |\cF_t\biggr].\ea\ee \et


%
\end{document}